\documentclass[12pt]{article}
\usepackage{amsmath,amssymb,amsbsy,amsfonts,amsthm,latexsym,
                        amsopn,amstext,amsxtra,euscript,amscd}

\newtheorem{lem}{Lemma}
\newtheorem{lemma}[lem]{Lemma}

\newtheorem{thm}{Theorem}
\newtheorem{theorem}[thm]{Theorem}
\newtheorem{cor}{Corollary}

 \DeclareMathOperator{\lcm}{lcm}

\def\mand{\qquad \mbox{and} \qquad}

\def\\{\cr}
\def\({\left(}
\def\){\right)}
\def\[{\left[}
\def\]{\right]}
\def\<{\langle}
\def\>{\rangle}
\def\fl#1{\left\lfloor#1\right\rfloor}
\def\rf#1{\left\lceil#1\right\rceil}

\def\cE{\mathcal E}

\def\cI{\mathcal I}

\def\cM{\mathcal M}
\def\cN{\mathcal N}

\def\cQ{\mathcal Q}
\def\cR{\mathcal R}

\def\cU{\mathcal U}
\def\cV{\mathcal V}

\def\e{{\mathbf{\,e}}}

\def\eps{\varepsilon}
\def\R{\mathbb{R}}

\def\T{\mathcal{T}}

\def\Z{\mathbb{Z}}

\def\tsigma{\widetilde{\sigma}}

\begin{document}

\title{Uniform Distribution of Fractional Parts Related to Pseudoprimes}
\author{
{\sc William D.~Banks} \\
{Department of Mathematics} \\
{University of Missouri} \\
{Columbia, MO 65211 USA} \\
{\tt bbanks@math.missouri.edu} \\
\and
{\sc Moubariz Z.~Garaev} \\
{Instituto de Matem{\'a}ticas}\\
{ Universidad Nacional Aut\'onoma de M{\'e}xico} \\
{C.P. 58089, Morelia, Michoac{\'a}n, M{\'e}xico} \\
{\tt garaev@matmor.unam.mx} \\
\and
{\sc Florian~Luca} \\
{Instituto de Matem{\'a}ticas}\\
{ Universidad Nacional Aut\'onoma de M{\'e}xico} \\
{C.P. 58089, Morelia, Michoac{\'a}n, M{\'e}xico} \\
{\tt fluca@matmor.unam.mx} \\
\and
{\sc Igor E.~Shparlinski} \\
{Department of Computing}\\
{Macquarie University} \\
{Sydney, NSW 2109, Australia} \\
{\tt igor@ics.mq.edu.au}}

\date{\today}

\pagenumbering{arabic}

\date{\today}
\maketitle

\newpage

\begin{abstract}
We estimate exponential sums with the Fermat-like quotients
$$
f_g(n) = \frac{g^{n-1} - 1}{n} \mand
h_g(n)=\frac{g^{n-1}-1}{P(n)},
$$
where $g$ and $n$ are positive integers, $n$ is composite, and
$P(n)$ is the largest prime factor of $n$. Clearly, both $f_g(n)$
and $h_g(n)$ are integers if $n$ is a Fermat pseudoprime to base
$g$, and if $n$ is a Carmichael number this is true for all $g$
coprime to $n$. Nevertheless, our bounds imply that the fractional
parts $\{f_g(n)\}$ and $\{h_g(n)\}$ are uniformly distributed, on
average over~$g$ for $f_g(n)$, and individually for $h_g(n)$. We
also obtain similar results with the functions ${\widetilde
f}_g(n) = gf_g(n)$ and  ${\widetilde h}_g(n) = gh_g(n)$.
\end{abstract}


\section{Introduction}

Throughout the paper, we use $P(n)$ to denote the largest prime
divisor of the integer $n\ge 2$, and we put $P(1)=1$.

For every integer $g\ge 1$, let $f_g(\cdot)$ and $h_g(\cdot)$ be
the arithmetic functions defined by
$$
f_g(n)=\frac{g^{n-1}-1}{n}\mand
h_g(n)=\frac{g^{n-1}-1}{P(n)}\qquad (n\ge 1).
$$
Clearly, $f_g(n)$ and $h_g(n)$ are \emph{integers} if $n$ is a
prime number and $n\nmid g$.  On the other hand, if $n$ takes only
\emph{composite} values, the problem of understanding the
distribution of the fractional parts of $f_g(n)$ and $h_g(n)$ is
rather involved. To approach this problem, we consider exponential
sums of the form:
\begin{eqnarray*}
S_g(a;N) &=& \sum_{\substack{n=1\\
n~{\mathrm{composite}}}}^N \e(a h_g(n)),\\
W(a;N) &=& \sum_{\substack{n=1\\
n~{\mathrm{composite}}}}^N~\sum_{\substack{g=1\\ \gcd(g,n)
=1}}^n\e(a f_g(n)),
\end{eqnarray*}
where the additive character $\e(\cdot)$ is defined (as usual) by
$\e(x)=\exp(2\pi ix)$ for all $x\in\R$, and $a\ne 0$ is an
integer.

We also consider the arithmetic functions
$$
{\widetilde f}_g(n) = \frac{g^{n} - g}{n} \mand {\widetilde
h}_g(n) =  \frac{g^{n} - g}{P(n)}\qquad (n\ge 1)
$$
and the corresponding exponential sums
\begin{eqnarray*}
{\widetilde S}_g(a;N) &=& \sum_{\substack{n=1\\
n~{\mathrm{composite}}}}^N \e(a {\widetilde h}_g(n)),\\
{\widetilde W}(a;N) &=& \sum_{\substack{n=1\\
n~{\mathrm{composite}}}}^N~\sum_{\substack{g=1\\ \gcd(g,n)
=1}}^n\e(a {\widetilde f}_g(n)).
\end{eqnarray*}

Clearly, ${\widetilde S}_g(a;N) = S_g(ag;N)$; the sums
${\widetilde W}(a;N)$, however, require an independent treatment.

Our results imply that the fractional parts $\{f_g(n)\}$,
$\{{\widetilde f}_g(n)\}$, $\{h_g(n)\}$ and
  $\{{\widetilde h}_g(n)\}$ are \emph{uniformly distributed}
over the interval $[0,1)$, on average over $g\in(\Z/n\Z)^*$ for
$f_g(n)$ and ${\widetilde f}_g(n)$, and individually (that is,
with $g>1$ fixed) for $h_g(n)$ and ${\widetilde h}_g(n)$. Of
course, one can either include or exclude the prime numbers in the
preceding statement since their contribution cannot change the
property of uniform distribution.

We remark that if $n$ is a \emph{Fermat pseudoprime to base $g$},
then both $f_g(n)$ and $h_g(n)$ are integers. If $n$ is a
\emph{Carmichael number}, then it is a Fermat pseudoprime to base
$g$ for every $g$ coprime to $n$, hence $f_g(n)$ and $h_g(n)$ are
integers for all such $g$. Since it is expected that there are
$$
C(N) = N^{1-(1+o(1))\log\log\log N/\log\log N}
$$
Carmichael numbers $n\le N$ (see~\cite{AGP,GrPom}), their
contribution to the sums $S_g(a;N)$ and $W(a;N)$ is substantial;
therefore, one cannot expect to obtain very strong bounds for
those sums. In particular, it is unlikely that one can obtain
upper bounds for  $S_g(a;N)$ and $W(a;N)$  of the form $O(N^{1
+\theta})$ and  $O(N^\theta)$, respectively, for any fixed
constant $\theta<1$. Indeed, using the \emph{Erd\H os-Tur{\'a}n
inequality}, which relates exponential sums to uniformity of
distribution, we show that the lower bound $S_g(a;N)\gg N/\log N$
holds for at least one integer $a$ in the range $1\le a\le\log N$;
thus, our upper bound for $S_g(a;N)$ (cf.~Theorem~\ref{thm:Sum S})
is rather tight. The same comments certainly apply to ${\widetilde
S}_g(a;N)$ and
  ${\widetilde W}(a;N)$ as well.

Problems of a similar flavor concerning the integrality and the
distribution of fractional parts of ratios formed with various
number theoretic functions have been treated previously
in~\cite{BGLS,BEPS,Luc,LS,Spiro1,Spiro2}. In part, our motivation
also stems from the results of~\cite{HeBr,HBK} on bounds for
exponential sums with \emph{Fermat quotients}.

It is perhaps surprising that, in order to establish our upper
bounds for $S_g(a;N)$ and $W(a;N)$, we need to apply tools from
very different and seemingly unrelated areas of number theory,
including several recent results. For instance, we not only apply
an asymptotic formula for the number of solutions to a symmetric
equation with an exponential function, which dates historically
back to 1962 (see the corollary to~\cite[Lemma~1,
Chapter~15]{Post}), but we also use very recent results on short
exponential sums from~\cite{Bour,BGK}. In the course of our
proofs, we also establish several new auxiliary results which may
be of independent interest; see, for example, Lemmas~\ref{lem:Bad
M-2} and~\ref{lem:AverageShortSums}.

In what follows, we use the Landau symbols $O$ and $o$, as well as
the Vinogradov symbols $\ll$ and $\gg$, with their usual meanings.
Any implied constants may depend, where obvious, on the parameter
$g$  but are absolute otherwise. We recall that the notations
$A\ll B$, $B\gg A$, and $A=O(B)$ are all equivalent, and $A=o(B)$
means that $A/B$ tends to zero. Throughout, we use the letters $p$
and $q$ exclusively to denote prime numbers, while $m$ and $n$
always denote positive integers. For a positive real number $x$ we
write $\log x$ for the maximum between the natural logarithm of
$x$ and $1$.

\section{Preliminary Results}

\subsection{Arithmetic Estimates}

Recall that a positive integer $n$ is said to be \emph{$y$-smooth}
if $P(n)\le y$.  For real numbers $x\ge y\ge 2$, let
$$
\Psi(x,y)=\#\{n\le x~:~P(n)\le y\}.
$$

\begin{lemma}
\label{lem:Smooth} Let $u=(\log x)/(\log y)$, where $x\ge y\ge 2$.
If $u\to\infty$ as $x\to\infty$, and $u\le y^{1/2}$, then the
following estimate holds:
$$
\Psi(x,y) = xu^{-u+o(u)}.
$$
\end{lemma}

For a proof of the Lemma~\ref{lem:Smooth}, we refer the reader
to~\cite[Section~III.5.4]{Tenen}; we remark that the condition
  $u\le y^{1/2}$ can be relaxed slightly, but the statement of
Lemma~\ref{lem:Smooth} is sufficient for our purposes.

For every positive integer $n$, let $\rho(n)$ denote the largest
squarefree divisor~$r$ of $n$ for which $\gcd(r,n/r)=1$; then
$s=n/\rho(n)$ is the largest \emph{powerful} divisor of $n$
(recall that a positive integer $m$ is said to be powerful if
$p^2\mid m$ for every prime $p$ that divides $m$).

We need the following statement, which is~\cite[Lemma~7]{ESW}:

\begin{lemma}
\label{lem:Bad M-1} Uniformly for $x\ge y\ge 1$, the bound
$\rho(n)>n/y$ holds for all $n\le x$ with at most $O(x/y^{1/2})$
exceptions.
\end{lemma}

For every positive integer $n$, let
$$
\gamma(n)=\prod_{p\,\mid\,n}\gcd(n-1, p-1).
$$
We note that this function also gives the cardinality of the set
of the so-called \emph{false witnesses} modulo $n$, that is, of
the set
$$
\{u \in \Z/n\Z~:~ u^{n-1} \equiv 1 \pmod n\},
$$
and has been studied in the literature (see~\cite{EP} and
references therein). The average value,  the normal order, and the
number of prime factors of $\gamma(n)$ are estimated in~\cite{EP};
however, these bounds do not seem to be enough for our purposes.

Our  next result   shows  for almost all \emph{composite} integers
$n$, the value of $\gamma(n)$ is very small. Although several
bounds on the number of composite integers $n\le x$ such that
$\gamma(n)>z$ can be extracted from~\cite{EP}, our estimate
appears to be new. More precisely,  \cite[Theorem~2.2]{EP} implies
such a bound for large values of $z$, and~\cite[Theorem~6.5]{EP}
treats the case of small values of $z$. In our applications,
however, we need a bound in the medium range. For our application,
it is convenient  to formulate this result in the following
two-parametric form:

\begin{lemma}
\label{lem:Bad M-2} Uniformly for  $x\ge y\ge 1$ and $\log\log\log
x=o(\log k)$, the number of composite integers $n\le x$ such that
$\gamma(n)>y^k$ is at most
$$
O\left(\frac{x\log\log x}{y}+\frac{x}{\exp((1+o(1))k\log
k)}\right).
$$
\end{lemma}

\begin{proof}
Let $\omega(m)$ be the number of distinct prime factors of the
$m$, and put $\cE_1=\{n\le x~:~\omega(n)\ge k\}$. If $n\in \cE_1$,
there exists a divisor $m\mid n$ with $\omega(m)=k$. For fixed
$m$, there are at most $x/m$ integers $n\in\cE_1$ such that
  $m\mid n$. Therefore, by unique factorization and the \emph{Stirling
formula} for $k!$, we see that
\begin{equation}
\label{eq:A1}
\begin{split}
\#\cE_1 & \le  x\sum_{\substack{m\le
x\\\omega(m)=k}}\frac{1}{m}\le\frac{x}{k!}\left(\,
\sum_{p^{\alpha}\le x}\frac{1}{p^{\alpha}}\right)^k=
\frac{x}{k!}\,(\log\log x+O(1))^k \\
& \le  x\left(\frac{e\log\log
x+O(1)}{k}\right)^k=x\exp\(-(1+o(1))k\log k)\),
\end{split}
\end{equation}
where the last estimate above uses the fact that
  $\log\log\log x=o(\log k)$.

Let $\varphi(\cdot)$ denote the \emph{Euler function}.
  We recall the estimate
\begin{equation}
\label{eq:Inverse Prime} \sum_{\substack{p\le t\\p\equiv 1\pmod
d}}\frac{1}{p} \ll \frac{\log\log t}{\varphi(d)},
\end{equation}
which holds uniformly for $1\le d\le t$ (see~\cite[Lemma~1]{Bas}
or the bound~(3.1) in~\cite{EGPS}). We also note that the bound
\begin{equation}
\label{eq:Inverse Euler}
\sum_{d>t}\frac{1}{d\varphi(d)}\ll\frac{1}{t}
\end{equation}
follows by partial summation from the asymptotic formula of
Landau~\cite{Lan}:
$$
\sum_{d\le t} \frac{1}{\varphi(d)} =
\frac{\zeta(2)\zeta(3)}{\zeta(6)}\left( \log t + \gamma -
\sum_{p}\frac{\log p}{p^2-p+1} \right ) + O\( \frac{\log t}{t} \),
$$
where $\zeta(s)$ is the \emph{Riemann zeta-function}, and $\gamma$
is the \emph{Euler-Mascheroni constant} (a more recent reference
is~\cite{Mont}).

Now let $\cE_2$ be the set of composite $n\le x$ for which there
exists $p\mid n$ with $d=\gcd(n-1,p-1)>y$. Write $n=pm$. Since
$n\equiv p\equiv 1\pmod d$, it follows that $m\equiv 1\pmod d$;
moreover, $m>1$ since $n$ is not prime. For each $p$ and $d$, we
have $1<m\le x/p$ and also $m \equiv 1\pmod d$, hence the number
of such $m$ is at most $x/pd$. Summing first over primes $p\equiv
1\pmod d$, then over all $d> y$, we derive from~\eqref{eq:Inverse
Prime} and~\eqref{eq:Inverse Euler} that
\begin{equation}
\label{eq:A2} \cE_2  \le \sum_{d> y}\sum_{\substack{p\le x\\
p\equiv 1\pmod d}}\frac{x}{pd} \ll x\sum_{d>y} \frac{\log\log
x}{d\varphi(d)}\ll \frac{x\log\log x}{y}.
\end{equation}

The result now follows from the estimates~\eqref{eq:A1}
and~\eqref{eq:A2} by observing that
$$
\gamma(n)=\prod_{p\,\mid\,n}\gcd(n-1,p-1)\le y^{\omega(n)}\le y^k
$$
if $n\le x$ is composite and not in the set $\cE_1\cup\cE_2$.
\end{proof}

By optimizing the choice of $y$ and $k$ for each given $z$, one
can reformulate Lemma~\ref{lem:Bad M-2} as the following more
concise (albeit weaker) statement:

\begin{cor}
Uniformly for $x\ge z\ge 1$ and $\log\log\log x=o(\log \log z)$,
the number of composite integers $n\le x$ such that $\gamma(n)>z$
does not exceed
$$
x\exp\(-{\sqrt{(0.5+o(1))\log z\log\log z}}\,\).
$$
\end{cor}

\begin{proof}
Choose $k$ such that $k^2\log k=\log z$, and put $y=z^{1/k}$.
Then, using our hypotheses on $x$ and $z$, we see that the
conditions of Lemma~\ref{lem:Bad M-2} are met, and the corollary
follows immediately.
\end{proof}

For a fixed base $g\ge 2$ and any prime $p\nmid g$, let $t_p$
denote the multiplicative order of $g$ modulo $p$. As usual, we
use $\tau(n)$ to denote the number of positive integer divisors of
$n$.

Let $\cQ$ be the set of primes $p$ satisfying the conditions
\begin{equation}
\label{eq:Set Q} \tau(p-1) \le (\log p)^2\mand t_p > p^{1/2} (\log
p)^{10},
\end{equation}
and let
\begin{equation}
\label{eq:Set R}
        \cR = \{p\ \text{prime} ~:~p\not\in\cQ\}.
\end{equation}

\begin{lemma}
\label{lem:Bad Primes} Uniformly for $x\ge 2$, the following bound
holds:
$$
\#\{p\le x~:~p\in \cR\}  \ll \frac{x}{(\log x)^2}.
$$
\end{lemma}

\begin{proof}
The result follows immediately from the \emph{Titchmarsh bound}:
$$
\sum_{p\le x}\tau(p-1) \ll x
$$
(see~\cite[Theorem~7.1, Chapter~5]{Prach})
and~\cite[Corollary~6]{IndlTim} (see also~\cite{ErdMur,Ford}).
\end{proof}

Finally, we need the following estimate:

\begin{lemma}
\label{lem:P(n) in interv} Let $A>0$ be fixed. Then, uniformly for
$x\ge y\ge 2$ and $\Delta>(\log y)^{-A}$, the following bound
holds:
$$
\#\{n\le x~:~y<P(n)\le y (1+ \Delta)\} \ll
\frac{x\log(1+\Delta)}{\log y},
$$
where the implied constant depends only on $A$.
\end{lemma}

\begin{proof} We can assume that $A$ is an integer (otherwise, replace it
with $\fl{A}$). We apply the following precise version of the
\emph{Mertens formula}:
\begin{equation}
\label{eq:strong_Mertens} \sum_{p\le t}\frac{1}{p}=\log\log
t+c_0+\frac{c_1}{\log t}+\cdots+ \frac{c_A}{(\log
t)^A}+O\(\frac{1}{(\log t)^{A+1}}\)
\end{equation}
for some constants $c_0,\ldots,c_A$, which follows by partial
summation from the \emph{Prime Number Theorem} (see, for example,
  \cite[Theorem~3.3, Chapter~3]{Prach}).
Applying~\eqref{eq:strong_Mertens} with $t=y$ and $t=y(1+\Delta)$,
and observing that for each prime $p$ in the interval
$\big(y,y(1+\Delta)\big]$, the number of integers $n\le x$ with
$P(n)=p$ does not exceed $x/p$, we obtain that
\begin{eqnarray*}
\lefteqn{\frac{1}{x}\cdot\#\{n\le x~:~y<P(m)\le y(1+\Delta)\}\le
\sum_{y<p\le y(1+\Delta)}\frac{1}{p}}\\
&&=\log\left(\log y+\log (1+\Delta)\right)-\log\log y\\
&&\qquad+\sum_{j=1}^A c_j\left(\frac{1}{(\log (y(1+\Delta)))^j}-
\frac{1}{(\log y)^j}\right)+O\left(\frac{1}{(\log y)^{A+1}}\right)\\
&&=\log\left(1+\frac{\log(1+\Delta)}{\log y}\right) +
O\left(\left|\frac{1}{\log (y(1+\Delta))}-\frac{1}{\log y}\right|+
\frac{1}{(\log y)^{A+1}}\right)\\
&&=\log\left(1+\frac{\log(1+\Delta)}{\log y}\right)+
O\left(\frac{\log(1+\Delta)}{(\log y)^2}+\frac{1}{(\log
y)^{A+1}}\right).
\end{eqnarray*}
If $\Delta$ is small, the first term above is approximately
$\Delta/\log y\gg(\log y)^{-(A+1)}$; hence, the error term never
dominates, and the result follows.
\end{proof}

\subsection{Estimates for Exponential Sums}

We begin with some well known and elementary results.

The following result, based on the \emph{Chinese Remainder
Theorem}, allows one to reduce exponential sums with polynomials
and with arbitrary denominators to exponential sums with prime
power denominators; this has been discussed, for example,
in~\cite[Problem 12.d, Chapter 3]{Vinog}:

\begin{lemma}
\label{lem:CRT} Let $n=n_1n_2$, where $n_1,n_2\ge 2$ are coprime,
and suppose that the integers $r_1,r_2$ satisfy:
\begin{eqnarray*}
r_1n_2 \equiv 1 \pmod {n_1} \mand r_2n_1 \equiv 1 \pmod {n_2}.
\end{eqnarray*}
Then, for any polynomial $F(X) \in \Z[X]$ with integer
coefficients, we have
$$
\sum_{\substack{g=0 \\ \gcd(g,n)=1}}^{n-1} \e\(F(g)/n\) =
\sum_{\substack{g_1=0 \\ \gcd(g_1,n_1)=1}}^{n_1-1}
\e\(r_1F(g_1)/n_1\) \sum_{\substack{g_2=0 \\
\gcd(g_2,n_2)=1}}^{n_2-1} \e\(r_2F(g_2)/n_2\).
$$
\end{lemma}

\begin{lemma}
\label{lem:Weil-CRT} For integers $a,n,k$ with $n,k\ge 1$, we have
$$
\left|\sum_{\substack{g=0\\ \gcd(g,n)=1}}^{n-1}
\e\(ag^k/n\)\right|\le  n d^{1/2} \gamma(n)\rho(n)^{-1/2},
$$
where  $d=\gcd(a, n)$.
\end{lemma}

\begin{proof} The proof is similar to that of~\cite[Lemma~4]{ESW}.
We recall the \emph{Weil bound}, which asserts that for every
integer $b$ and prime $p\nmid b$, the inequality
$$\left|\sum_{g=1}^{p-1} \e\(bg^k/p\)\right|\le
\gcd(k,p-1)  p^{1/2}
$$
holds (see, for example, \cite[Theorem~5.41]{LN}).

Let $\rho(n)=p_1\ldots p_\nu$ be the factorization of $\rho(n)$ as
a product of (distinct) primes, and put $s = n/\rho(n)$. Then, by
Lemma~\ref{lem:CRT}, we have
$$
\sum_{\substack{g=0\\ \gcd(g,n)=1}}^{n-1} \e\(ag^k/n\) =
\prod_{j=1}^\nu \left(\,\sum_{g_j=1}^{p_j-1} \e_{p_j}\(a
b_jg_j^k/p_j\)\right) \left(\,\sum_{\substack{h=0\\
\gcd(h,s)=1}}^{s-1} \e\(ac h^k/s\)\right)
$$
for some integers $b_1, \ldots, b_\nu$ and $c$ such that
$\gcd(b_j,p_j)=1$ for $j =1,\ldots,\nu$ and $\gcd(c,s)=1$. For
each $j$ such that $p_j\mid a$, the sum over $g_j$ is equal to
$p_j-1$. We estimate the sum over $h$ trivially as $s$. Therefore,
$$
\left|\sum_{\substack{g=0\\ \gcd(g,n)=1}}^{n-1}
\e\(ag^k/n\)\right| \le s \prod_{\substack{j=1\\ p_j\,\nmid\,
a}}^\nu
\left(\gcd(k, (p_j-1) p_j^{1/2}\right) \prod_{\substack{j=1\\
p_j\,\mid\,a}}^\nu p_j,
$$
and the result follows.
\end{proof}

The next result appears in~\cite{Bour}; it can also be deduced
from~\cite[Theorem~5]{BGK} in an even more explicit form:

\begin{lemma}
\label{lem:ShortSums} For every $\delta > 0$, there exists
$\eta>0$, such that if
$$
p^\delta\le M\le t_p,
$$
then for every integer $a$ not divisible by $p$, the following
bound holds:
$$
\left| \sum_{m \le M} \e(a g^m/p) \right| \le Mp^{-\eta}.
$$
\end{lemma}

The following bound on short exponential sums with an exponential
function appears to be new and may be of independent interest. To
prove this bound, we use the well known method of estimating
double exponential sums via the number to solutions of certain
symmetric systems of equations, which can be found
in~\cite{GLS,Kar1,Kar2,Kar3,KonShp,Kor1,Kor2} and in many other
places (see, for example, \cite{Kon}). In fact, although the
result is conveniently summarized in~\cite[Lemma~4]{Kon}, no proof
is given there. Here, we supply a proof for the sake of
completeness.

\begin{lemma}
\label{lem:AverageShortSums} For a real number $V\ge 2$ and
positive integers $M,k,\ell$ satisfying the inequalities
$$
2^k k!\,\pi(V) \le M^{k+1}, \mand 2^{\ell} \ell!\,\pi(V) \le
M^{(\ell +1)/2},
$$
the following bound holds:
$$
\sum_{\substack{p \le V\\p\,\nmid\,ag}} ~\max_{L \le M} \left|
\sum_{m=1}^L \e(ag^m/p)\right| \ll \pi(V) M \(\frac{V^{1/2}
M^{3/4}}{\pi(V)}\)^{1/k\ell} ,
$$
where  the implied constant depends only on $g$.
\end{lemma}

\begin{proof} For each prime $p \le V$ such that $p\nmid ag$,
let $L_p$ denote the smallest positive integer such that
$$
\max_{L \le M} \left| \sum_{m=1}^L \e(ag^m/p)\right| =  \left|
\sum_{m=1}^{L_p} \e(ag^m/p)\right|.
$$
Put $H = \fl {M^{1/2}}$; then,
\begin{equation}
\label{eq:Sum W} \sum_{\substack{p \le V\\p\,\nmid\,ag}} \left|
\sum_{m=1}^{L_p} \e(ag^m/p)\right| = \frac{W}{H}  + O(\pi(V) H),
\end{equation}
where
$$
W = \sum_{\substack{p \le V\\p\,\nmid\,ag}} \sum_{h=1}^H \left|
\sum_{m=1}^{L_p} \e(ag^{m+h}/p)\right|.
$$
By the \emph{H{\"o}lder inequality}, it follows that
\begin{eqnarray*}
W^k & \le & \pi(V)^{k-1} H^{k-1} \sum_{\substack{p \le
V\\p\,\nmid\,ag}} \sum_{h=1}^H \left| \sum_{m=1}^{L_p}
\e(ag^{m+h}/p)\right|^k\\
& = & \pi(V)^{k-1} H^{k-1} \sum_{\substack{p \le V\\p\,\nmid\,ag}}
\sum_{h=1}^H \vartheta_{p,h} \( \sum_{m=1}^{L_p}
\e(ag^{m+h}/p)\)^k
\end{eqnarray*}
for some complex numbers $\vartheta_{p,h}$ of absolute value $1$.

Now, let $R_{p,s}(K, \lambda)$ denote the number of solutions of
the congruence
$$
\sum_{i=1}^s   g^{r_i} \equiv \lambda \pmod p  \qquad (1 \le r_1,
\ldots, r_s \le K).
$$
Then
$$
               \(\,\sum_{m=1}^{L_p}
\e(ag^{m+h}/p)\)^k  = \sum_{\lambda =0}^{p-1} R_{p,k}(L_p,
\lambda) \e(a \lambda g^h/p).
$$
Therefore, after  changing the order of summation, we derive that
$$
W^k \le \pi(V)^{k-1} H^{k-1} \sum_{\substack{p \le
V\\p\,\nmid\,ag}} \sum_{\lambda =0}^{p-1} R_{p,k}(L_p, \lambda)
\sum_{h=1}^H \vartheta_{p,h} \e(a \lambda g^h/p).
$$
Writing
$$
R_{p,k}(L_p, \lambda) = \(R_{p,k}(L_p, \lambda)^2\)^{1/2\ell}
R_{p,k}(L_p, \lambda)^{(\ell -1)/\ell}
$$
and using the H{\"o}lder inequality for a sum of products of three
terms, we have
\begin{eqnarray*}
W^{2k\ell} & \le & \pi(V)^{2\ell(k-1)} H^{2\ell(k-1)}
\sum_{\substack{p \le V\\p\,\nmid\,ag}} \sum_{\lambda =0}^{p-1}
R_{p,k}(L_p, \lambda)^2\\
& & \qquad \qquad\times\quad \(\,\sum_{\substack{p \le
V\\p\,\nmid\,ag}} \sum_{\lambda =0}^{p-1}
R_{p,k}(L_p, \lambda)\)^{2 \ell -2}\\
& & \qquad \qquad  \times\quad \sum_{\substack{p \le
V\\p\,\nmid\,ag}} \sum_{\lambda =0}^{p-1} \left|\sum_{h=1}^H
\vartheta_{p,h} \e(a \lambda g^h/p)\right|^{2\ell}.
\end{eqnarray*}
Clearly,
$$
\sum_{\lambda =0}^{p-1} R_{p,k}(L_p, \lambda) = L_p^k \le M^k,
$$
and
$$
\sum_{\lambda =0}^{p-1} R_{p,k}(L_p, \lambda)^2 = T_{p,k}(L_p),
$$
where $T_{p,s}(K)$ denotes the number of solutions of the
congruence
$$
\sum_{i=1}^{2s}  (-1)^i g^{r_i} \equiv 0 \pmod p  \qquad (1 \le
r_1, \ldots, r_s \le K).
$$
Thus,
\begin{eqnarray*}
W^{2k\ell} & \le & \pi(V)^{2\ell(k-1) +  2\ell-2} H^{2\ell(k-1)}
M^{2k(\ell-1)}\sum_{\substack{p \le V\\p\,\nmid\,ag}} T_{p,k}(L_p)\\
&& \qquad \times \sum_{\substack{p \le V\\p\,\nmid\,ag}}
\sum_{\lambda =0}^{p-1} \left|\sum_{h=1}^H \vartheta_{p,h} \e(a
\lambda g^h/p)\right|^{2\ell}.
\end{eqnarray*}
Furthermore,
\begin{eqnarray*}
\lefteqn{\sum_{\lambda =0}^{p-1} \left|\sum_{h=1}^H
\vartheta_{p,h} \e(a \lambda g^h/p)\right|^{2\ell}= \sum_{h_1,
\ldots, h_{2\ell} =1}^H \prod_{i=1}^{2\ell} \vartheta_{p,h_i}
\sum_{\lambda =0}^{p-1} \e\(\frac{\lambda}{p}
\sum_{i=1}^{2\ell}  (-1)^i g^{h_i}\) }\\
&  & \qquad\qquad\qquad\qquad \le   \sum_{h_1, \ldots, h_{2\ell}
=1}^H \left| \sum_{\lambda =0}^{p-1} \e\(\frac{\lambda}{p}
\sum_{i=1}^{2\ell}  (-1)^i g^{h_i}\)  \right| =  p\,T_{p,\ell}(H).
\end{eqnarray*}
Hence,
$$
W^{2k\ell} \le \pi(V)^{2k \ell -2} H^{2\ell(k-1)}
M^{2k(\ell-1)}\sum_{\substack{p \le V\\p\,\nmid\,ag}} T_{p,k}(L_p)
\sum_{\substack{p \le V\\p\,\nmid\,ag}} p\, T_{p,\ell}(H).
$$
We remark that
$$
\sum_{\substack{p \le V\\p\,\nmid\,ag}} T_{p,k}(L_p)  \le \sum_{p
\le V } T_{p,k}(M)
$$
is equal to the number of primes $p \le V$ which divide all
possible expressions of the form
$$
\sum_{i=1}^{2k}  (-1)^i g^{m_i} \qquad (1 \le m_1, \ldots, m_{2k}
\le M).
$$
Clearly, any nonzero sum above has at most $\log (2k g^M)/\log 2$
prime divisors.  Also, by the corollary
  to~\cite[Lemma~1, Chapter~15]{Post}, there are at most $2^k k!
M^k$ such sums which vanish (see also~\cite{BNC} for a survey of
recent results in this direction). For these ones, we estimate the
number of prime divisors trivially as $\pi(V)$. Thus, using the
inequality $2^k k!\,\pi(V)\le M^{k+1}$, we deduce that
$$
\sum_{p \le V } T_{p,k}(M) \ll M^{2k+1} + M^{2k} \log k + 2^k
k!M^k \pi(V) \ll M^{2k+1} .
$$
Similarly,
$$
\sum_{\substack{p \le V\\p\,\nmid\,ag}} p T_{p,\ell}(H) \le V
\sum_{p \le V } T_{p,\ell}(H) \ll V H^{2\ell+1}.
$$
Consequently,
$$
W^{2k\ell} \ll \pi(V)^{2k \ell -2} V H^{2k \ell  + 1} M^{2k\ell
+1}.
$$
Substituting this estimate into~\eqref{eq:Sum W}, we obtain that
$$
\sum_{\substack{p \le V\\p\,\nmid\,ag}}  \left| \sum_{m=1}^{L_p}
\e(ag^m/p)\right| \ll \pi(V)^{1-1/k\ell} V^{1/2k\ell} M^{3/4k\ell}
+ \pi(V) M^{1/2}.
$$
It now remains only to observe that, since $2^k k!\,\pi(V) \le
M^{k+1}$, the last term never dominates.
\end{proof}

It is important to remark that the implied constant in the bound
of Lemma~\ref{lem:AverageShortSums} depends on $g$ but not on the
parameters $k,\ell$ (nor on $a,M,V$). In particular, in our
applications we can choose $k$ and $\ell$ to be growing functions
of $M$ and $V$. Of course, we use Lemma~\ref{lem:AverageShortSums}
only to deal with the case that $M$ is suitably small with respect
to $V$, and in the remaining range, we apply
Lemma~\ref{lem:ShortSums}.

We also need the following bound, which is a special case of the
more general results of~\cite{GarShp}. We recall that the set
$\cQ$ is defined by~\eqref{eq:Set Q}.

\begin{lemma}
\label{lem:LargeSieve} For any real number $U$, any positive
integer $M$, and any subset $\cM \subseteq \{1, \ldots, M\}$ of
cardinality $\# \cM = T$, we have the uniform bound:
$$
\sum_{\substack{p\in \cQ\\ U \le p \le 2U}} \max_{(a,
p)=1}\left|\sum_{m\in\cM} \e\(a g^{m}/p\)\right|^2\ll UT (M (\log
U)^{-20}  +    U) (\log U)^3.
$$
\end{lemma}

\section{Single Exponential Sums with $h_g(n)$}

\begin{theorem}
\label{thm:Sum S} Fix $g>1$ and $\eps>0$. Then for every integer
$a$ such that $ \log |a|\le \exp\( (\log N)^{1-\varepsilon}\)$,
the inequality
$$
S_g(a;N)\ll \frac{  N}{ \sqrt{\log  N}}
$$
holds, where the implied constant depends only on $g$ and $\eps$.
\end{theorem}

\begin{proof}
We may assume that $\varepsilon<1/2$. Put $Q =\exp\( 2(\log
N)^{1-\varepsilon}\)$, and let $\cE_1$ denote the set of
$Q$-smooth integers $n\le N$. Then, applying
Lemma~\ref{lem:Smooth} with $ u = 0.5 (\log N)^{\varepsilon}$, we
obtain the bound
\begin{equation}
\label{eq:E1}
\begin{split}
\# \cE_1 &=\Psi(N,Q)= N u^{-u+o(u)}\\
&=N\exp\(-\(0.5\,\varepsilon+o(1)\)(\log N)^{\varepsilon}\log\log
N\).
\end{split}
\end{equation}

Next, let $\cE_2$ be the set of the integers $n\le N$,
$n\not\in\cE_1$, such that $P(n)\mid ag$. We have
        \begin{equation}
\label{eq:E2} \# \cE_2 \le \sum_{\substack{p > Q\\p\,\mid\,ag}}
\frac{N}{p} \ll
        \frac{N}{Q} \sum_{p\,\mid\,ag} 1 \ll \frac{N}{Q} \log |a|
\le N \exp\( -(\log N)^{1-\varepsilon}\).
        \end{equation}

Let $\cE_3$ be the set of the positive integers $n\le N$ not in
$\cE_1$  such that $P(n) \in \cR$ where the set $\cR$ is defined
by~\eqref{eq:Set R}. We have
\begin{equation}
\label{eq:E3} \# \cE_3 \le \sum_{\substack{Q < p \le N\\ p \in
\cR}} \sum_{\substack{n \le N\\ P(n) = p}} 1 \le N
\sum_{\substack{Q < p \le N\\ p \in \cR}} \frac{1}{p}.
\end{equation}
By Lemma~\ref{lem:Bad Primes} and partial summation, we obtain
that
$$
\# \cE_3  \ll \frac{N}{\log Q}\le \frac{N}{\sqrt{\log N}}.
$$

Let us now denote
$$
X = N^{1/2}(\log N)^{-5}, \qquad Y = N^{3/4}\qquad {\text{\rm and
}} \qquad Z = N \exp\(- \sqrt{\log N}\,\).
$$

Let $\cE_4$ be the set of the positive integers $n\le N$ such that
either
$$
X<P(n) \le N^{1/2},
$$
or
$$
Z<P(n) \le N.
$$
By Lemma~\ref{lem:P(n) in interv}, it follows that
\begin{equation}
\label{eq:E4} \cE_4 \ll \frac{N}{\sqrt{\log N}}.
\end{equation}

Let $\cN$ be the set of integers $n\le N$ such that
$n\not\in\cE_1\cup\cE_2\cup\cE_3\cup\cE_4$. Then, from the
estimates \eqref{eq:E1}, \eqref{eq:E2}, \eqref{eq:E3} and
\eqref{eq:E4}, we conclude that
$$
S_g(a;N) = \sum_{n =1}^N \e(a h_g(n)) + O\(\frac{N}{\log N}\) =
\sum_{n \in \cN} \e(a h_g(n)) + O\(\frac{N}{\sqrt{\log N}}\).
$$
Note that the error term in the middle expression comes from prime
values of $n\le N$, which are not included in the sum $S_g(a;N)$.

Every $n \in \cN$ has a unique representation of the form $n=pm$,
with a prime $p \ge Q$ and an integer $m\le N/p$ such that $P(m)
\le p$. Also, remarking that for $p > N^{1/2}$ the condition $P(m)
\le p$ is automatically satisfied, we see that
$$
\sum_{n \in \cN} \e(a h_g(n)) = W_1 + W_2 + W_3,
$$
where, since $g^{pm} \equiv g^m \pmod p$, we have
\begin{eqnarray*}
|W_1 | & =& \left| \sum_{\substack{Q < p \le X\\ p \in \cQ}}
\sum_{\substack{m \le N/p\\ P(m) \le p}} \e(a h_g(pm))\right| \le
\sum_{\substack{Q < p \le X\\ p \in \cQ}} \left| \sum_{\substack{m
\le N/p\\ P(m) \le p}}
\e(ag^{m-1}/p)\right|,\\
|W_2| & =& \left| \sum_{\substack{N^{1/2} < p \le Y\\ p \in \cQ}}
\sum_{m \le N/p} \e(a h_g(pm))\right| \le \sum_{\substack{N^{1/2}
< p \le Y\\ p \in \cQ}} \left|\sum_{m \le N/p}
\e(ag^{m-1}/p)\right|,\\
|W_3| & =& \left| \sum_{\substack{Y < p \le Z\\ p \in \cQ}}
\sum_{m \le N/p} \e(a h_g(pm))\right| \le \sum_{\substack{Y < p
\le Z\\ p \in \cQ}} \left|\sum_{m \le N/p} \e(ag^{m-1}/p)\right| .
\end{eqnarray*}

To estimate $|W_1|$, put $\Delta = 1/\log N$ and consider the
sequence of real numbers:
$$
U_j = \min\{Q (1 + \Delta)^j, X\} \qquad (0\le j\le J),
$$
where
\begin{equation}
\label{eq:J bound} J=\rf{\frac{\log (X/Q)}{\log (1 + \Delta)}} \ll
\Delta^{-1} \log N=(\log N)^{2}.
\end{equation}
We denote the set of primes $p \in \cQ$ in the half-open interval
$(U_j, U_{j+1}]$ by $\cU_j$, $j =0, \ldots, J-1$. Note that since
$$
\Delta=(\log N)^{-1}\ge (\log Q)^{-2}\ge (\log U_j)^{-2},
$$
we can apply Lemma~\ref{lem:P(n) in interv} with $A=2$ in what
follows. {From} the above, we infer that
\begin{equation}
\label{eq:split-1} |W_1| \le \sum_{j=0}^{J-1} | \sigma_j | ,
\end{equation}
where
$$
\sigma_j=\sum_{p\in\cU_j} \sum_{\substack{m\le N/p\\
P(m) \le p}}\e(ag^{m-1}/p)\qquad (0\le j\le J-1).
$$
We have
\begin{eqnarray*}
\sigma_j & =&  \sum_{p \in \cU_j} \(\sum_{\substack{m \le
N/U_j\\P(m) \le p}}
\e(ag^{m-1}/p)  + O\big(|N/p - N/U_j|\big)\)\\
& =&  \sum_{p \in \cU_j} \(\sum_{\substack{m \le N/U_j\\ P(m) \le
p}} \e(ag^{m-1}/p)  + O(N\Delta/p)\).
\end{eqnarray*}
Applying Lemma~\ref{lem:P(n) in interv} and using the fact that
$\log(1+\Delta)\le \Delta$, we obtain that
\begin{eqnarray*}
\sigma_j & =&  \sum_{p \in \cU_j} \(\sum_{\substack{m \le N/U_j\\
P(m) \le U_j}}
\e(ag^{m-1}/p)  + O\big(N\Delta/p + N\Delta/(U_j \log U_j)\big)\)\\
               & =&  \tsigma_j  + O\(N\Delta \sum_{p \in \cU_j} 1/p\),
\end{eqnarray*}
where
$$
\tsigma_j = \sum_{p \in \cU_j} \sum_{\substack{m \le N/U_j\\
P(m) \le U_j}}  \e(ag^{m-1}/p)\qquad (0\le j\le J-1).
$$
Thus, from~\eqref{eq:split-1}, we have
\begin{equation}
\label{eq:split-2} | W_1| \le \sum_{j=0}^{J-1} | \tsigma_j | + O\(
N\Delta \sum_{p \le N} 1/p\) = \sum_{j=0}^{J-1} | \tsigma_j | +
O\( \frac{N\log \log N}{\log N}\).
\end{equation}
Using the trivial bound  $\#\cU_j\le\Delta U_j$ $\big(\,$in fact,
the stronger bound
$$
\#\cU_j\ll \Delta U_j/\log U_j\le \Delta U_j/\log Q
$$
also holds (see~\cite{MV}, for example), but this does not lead to
an improvement in the final bound for $S_g(a;N)\big)$ and the
\emph{Cauchy inequality}, we derive that
$$
\tsigma_j^2 \le\Delta U_j \sum_{p \in \cU_j} \left|
\sum_{\substack{m \le N/U_j\\ P(m) \le U_j}}
\e(ag^{m-1}/p)\right|^2.
$$
Applying Lemma~\ref{lem:LargeSieve} and estimating the number of
$m \le N/U_j$ such that $P(m) \le U_j$ trivially as $N/U_j$, we
see that
\begin{eqnarray*}
|\tsigma_j|^2 &\ll& \Delta N U_j (N U_j^{-1} (\log U_j)^{-20}  +
U_j) (\log
U_j)^3\\
&  = & \Delta  N^2  (\log U_j)^{-17}  + \Delta N U_j^2 (\log U_j)^3 \\
&\le& \Delta  N^2  (\log Q)^{-17} + \Delta N X^2 (\log N)^3\le
2N^2 (\log N)^{-8}.
\end{eqnarray*}
Therefore, from~\eqref{eq:J bound} and~\eqref{eq:split-2} it
follows that
$$
|W_1| \ll\frac{N\log\log N}{\log N}.
$$

To estimate $W_2$, we simply apply Lemma~\ref{lem:ShortSums} with
$\delta = 1/6$ to each sum over $m$, getting
$$
\sum_{m \le N/p} \e(ag^{m-1}/p) \ll \frac{N}{p}  p^{-\eta}
$$
with some absolute constant $\eta > 0$. Here, recall that $t_p \ge
p^{1/2}$ for every prime $p  \in \cQ$; hence, the above bound
follows from Lemma~\ref{lem:ShortSums} regardless of whether
$t_p\ge N/p$ or not.  Consequently,
$$
|W_2| \ll \sum_{ N^{1/2}<p \le N} \frac{N}{p}  p^{-\eta} \le
N^{1-\eta/2} \sum_{ N^{1/2}<p \le N} \frac{1}{p} \ll N^{1-\eta/2}
\log \log N.
$$

To estimate $W_3$, consider the sequence of real numbers:
$$
V_i=\max\{ Y , e^{-i} Z\}\qquad (0\le i\le I),
$$
where $I=\rf{\log (Z/Y)}$. We denote the set of primes $p \in \cQ$
in the half-open interval $(V_{i+1}, V_i]$ by $\cV_i$, $i =0,
\ldots, I-1$. Then
\begin{equation}
\label{eq:split-3} |W_3| \le \sum_{i=0}^{I-1} | \varSigma_i | ,
\end{equation}
where
$$
\varSigma_i = \sum_{p \in \cV_i}\sum_{m \le N/p}\e(ag^{m-1}/p).
$$
For each $i =0, \ldots, I-1$, we apply
Lemma~\ref{lem:AverageShortSums} with the parameter choices
$$
k =   \ell = \rf{\frac{4 \log N}{i + \sqrt{\log   N}}}, \qquad V =
V_{i+1} \qquad {\text{\rm and}} \qquad M = \rf{N/V_i}.
$$
In particular,
$$
M \ge  \exp\(i - 1 + \sqrt{\log   N} \),
$$
and also
$$
\frac{N}{\log N}\ll \pi(V)M\ll\frac{N}{\log N}.
$$
Since, for sufficiently large $N$, the inequality
$$
\frac{M^{(\ell+1)/2}}{2^{\ell}\ell!}\ge\frac{M^{k/2}}{2^kk!} \ge
\(\frac{M^{1/2}}{2k}\)^k \ge M^{k/3}\ge e^{(i+\sqrt{\log
N}\,)k/4}\ge N
$$
holds, one easily verifies  that the conditions of
Lemma~\ref{lem:AverageShortSums} are satisfied if $N$ is large
enough. Since $V>N^{3/4}$ and $M<N^{1/4},$ we have
$$
M^{3/4}V^{-1/2}\log V\ll N^{-3/16}\log N\ll N^{-1/6}.
$$
Thus, an application of Lemma~\ref{lem:AverageShortSums} yields
the bound
\begin{eqnarray*}
\left|\varSigma_i\right| &\ll &\frac{N}{\log
N}\(N^{-1/6}\)^{1/k\ell}= \frac{N}{\log N}\exp\(-\frac{1}{150}
\(i+\sqrt{\log N}\,\)^2/\log N\)
\\ & \le &\frac{N}{\log N}\exp\(-\frac{i^2}{150\log N}\).
\end{eqnarray*}
From~\eqref{eq:split-3}, we now derive that
$$
|W_3| \le \frac{N}{\log N} \sum_{i=0}^{\infty}e^{-i^2/150\log
N}\ll \frac{N}{\log N} \int_{0}^{\infty}e^{-t^2/150\log N}d t \ll
\frac{N}{  (\log N)^{1/2}},
$$
and the proof is complete.
\end{proof}

Next, we obtain a lower bound which shows that the upper bound of
Theorem~\ref{thm:Sum S} is quite tight.

\begin{theorem}
\label{thm:Sum S Lower} Let $g>1$ be a fixed integer base. Then
the inequality
$$
\max_{1 \le a \le \log N} |S_g(a;N)|\gg \frac{  N}{  \log  N}
$$
  holds, where the implied constant depends only on $g$.
\end{theorem}

\begin{proof}
   Let $\T$ be the set of positive integers $n\le N$ which
can be expressed in the form $n=mp$, where the prime $p$ and
integer $m$ satisfy the inequalities
$$
m\le \frac{\log N}{6\log g}, \qquad N^{2/3}<p\le N/m.
$$
Clearly, for each $m$ there are $(1+o(1))N/(m\log N)$ primes $p$
such that $n=mp$ lies in $\T$, and the pair $(m,p)$ is uniquely
determined by $n$. Therefore,
$$
\# \T\gg \sum_{m\le (\log N)/(6\log g)}\frac{N}{m\log N}\gg
\frac{N\log\log N}{\log N}.
$$
Next, observe that for every $n\in \T$,
$$
\{h_g(n)\}=\left\{\frac{g^{mp-1}-1}{p}\right\}=\left\{\frac{g^{m-1}-1}{p}
\right\} <\frac{N^{1/6}}{N^{2/3}} = N^{-1/2}.
$$
Thus, the numbers $\{h_g(n)\}$ with $n\in \T$ all lie in the
interval $[0, N^{-1/2}).$

On the other hand, by the \emph{Erd\H os-Tur{\'a}n inequality}
(see~\cite[Theorem~1.21, Section~1.2.2]{DrTi}
or~\cite[Theorem~2.5, Section~2.2]{KuNi}) for the number of points
$A(\gamma)$ in an interval $[0,\gamma) \subseteq [0,1)$,
$$
\max_{0\le \gamma \le 1} |A(\gamma) - \gamma N| \ll
\frac{N}{H}+\sum_{a=1}^{H}\frac{1}{a}S_g(a, N)\qquad(H\ge 1).
$$
Therefore, applying this inequality with $\gamma = N^{-1/2}$, we
derive
$$
\frac{N\log\log N}{\log N}\ll \#\T\ll N^{1/2}+
\frac{N}{H}+\sum_{a=1}^{H}\frac{1}{a}S_g(a, N).
$$
Hence, by taking $H= \fl{\log N}$, and assuming that $N$ is large
enough, we obtain the stated result.
\end{proof}

It is easy to see that choosing a smaller value of $H$, one can
obtain the lower bound of Theorem~\ref{thm:Sum S Lower} over the
smaller range $1 \le a \le c(g)\log N/\log\log N$ for some
constant $c(g)>0$ depending only on $g$.

\section{Double Exponential  Sums with $f_g(n)$}

\begin{theorem}
\label{thm:Sum D} For any integer $a$ such that $\log |a|
=o\(\sqrt{\log N\log\log N}\,\)$, the following inequality holds:
$$
W(a;N)\le N^2 \exp\(-(0.5+o(1))\sqrt {\log N\log\log N}\,\).
$$
\end{theorem}

\begin{proof}
Let $N$ be sufficiently large, and suppose that $k$ (a positive
integer parameter that depends only on $N$) is such that
$\log\log\log N=o(\log k)$. Put $y=\exp(k\log k)$, and let $\cE$
be the set of composite integers $n\le N$ such that either
$\rho(n)\le n/y^2$ or $\gamma(n)>y^k$. By Lemmas~\ref{lem:Bad M-1}
and~\ref{lem:Bad M-2}, it follows that
$$
|W(a;N)|\le \sum_{\substack{n\le N,~n\not\in {\cE}\\
n~{\mathrm{composite}}}}\left|\sum_{\substack{g=1\\
\gcd(g,n)=1}}^n\e(af_g(n))\right|+
O\left(\frac{N^2}{\exp((1+o(1))k\log k)}\right).
$$
In  $n\not\in \cE$, then $\rho(n)>n/y^2$ and $\gamma(n)\le y^k$;
hence, by Lemma~\ref{lem:Weil-CRT}, we see that
\begin{eqnarray*}
|W(a;N)| & \ll & |a|y^{k+1}N^{3/2}+\frac{N^2}{\exp\((1+o(1))k\log k\)}\\
& = & |a|N^{3/2}\exp\(k(k+1)\log
k\)+\frac{N^2}{\exp\((1+o(1))k\log k\)}.
\end{eqnarray*}
Choosing $k$ such that $k(k+2)\log k=(0.5+o(1))\log N$ (to balance
the two terms above), we obtain the stated estimate.
\end{proof}

\section{Double Exponential  Sums with ${\widetilde f}_g(n)$}

\begin{theorem}
\label{thm: D'} For any nonzero integer $a$ with $|a| <
(\log\log\log N)^3$ the bound
$$
{\widetilde W}(a;N)\ll \frac{N^2\log\log\log\log N}{\log\log\log
N}
$$
holds as $N\to \infty$.
\end{theorem}

\begin{proof}
Let $\lambda(\cdot)$ denote the \emph{Carmichael function}. We
recall that if
$$
n=\prod_{\nu=1}^{s} p_\nu^{\alpha_\nu}
$$
is the prime factorization of $n$, then
$$
\lambda(n)=\lcm\left[\lambda(p_1^{\alpha_1}),
\ldots,\lambda(p_s^{\alpha_s}) \right],
$$
where $\lambda(p^{\alpha})=p^{\alpha-1}(p-1)$ for a prime power
except when $p=2$ and $\alpha\ge 3$, in which case
$\lambda(2^{\alpha})=2^{\alpha-2}$.

Put
$$
y=(\log\log\log N)^2\qquad {\text{\rm and}}\qquad z=\frac{\log\log
N}{(\log\log \log N)^2},
$$
and let $\cI$ be the interval $[y,z]$.

The proof of~\cite[Lemma~2]{LuPo} shows that if $\cE_1$ is the set
of integers $n\le N$ for which there exists a prime number
  $q\in \cI$ such that $q\nmid\lambda(n)$, then
\begin{equation}
\label{eq:A11} \#\cE_1\ll \frac{N}{\log\log N}.
\end{equation}
Let $\cE_2$ be the set of $n\le N$ such that $q^2\mid n$ for some
prime $q>y$. Then
\begin{equation}
\label{eq:A21} \#\cE_2\le \sum_{q\ge y}\frac{N}{q^2}\ll
\frac{N}{y}\ll \frac{N}{(\log\log\log N)^2}.
\end{equation}
Let $\cE_3$ be the set of $n\le N$ such that $n$ is not divisible
by any prime in $\cI$. By the inclusion-exclusion principle, we
have
\begin{equation}
\label{eq:A31}
\begin{split}
       \#\cE_3 &=  N\prod_{y\le q\le
z}\left(1-\frac{1}{q}\right) + O(2^z) \ll N\,\frac{\log y}{\log z} + 2^z\\
& \ll \frac{N\log\log\log\log N}{\log\log\log N}.
\end{split}
\end{equation}

Finally, let $\cN$ be the set of integers $n\le N$ such that
$n\not\in\cE_1\cup\cE_2\cup\cE_3$. Thus, from~\eqref{eq:A11},
\eqref{eq:A21} and~\eqref{eq:A31}, we deduce that
\begin{eqnarray}
\label{eq:W and sigma} {\widetilde
W}(a;N)=\sigma+O\left(\frac{N^2\log\log\log\log N}{\log\log\log
N}\right),
\end{eqnarray}
where
$$
\sigma=\sum_{n\in \cN}\sum_{\substack{g=1\\
\gcd(g,n) =1}}^n\e(a {\widetilde f}_g(n)).
$$
To handle this sum, write $d_n =\gcd (n, \lambda(n))$, and put
$s_n = \lambda(n)/d_n$. Then
\begin{eqnarray*}
\sigma & = & \sum_{n\in \cN} \sum_{\substack{g=1\\
\gcd(g,n) =1}}^n \e(a(g^n - g)/n)\\ & = & \sum_{n\in
\cN}\frac{1}{\varphi(n)}\sum_{\substack{1\le
h\le n\\ \gcd(h,n)=1}}\sum_{\substack{g=1\\
\gcd(g,n) =1}}^n
\e(a((gh^{s_n})^n - g h^{s_n})/n)\nonumber\\
& = & \sum_{n\in \cN}
       \frac{1}{\varphi(n)} \sum_{\substack{g=1\\
\gcd(g,n) =1}}^n  \sum_{\substack{1\le h\le n\\ \gcd(h,n)=1}}
\e(a(g^n - g h^{s_n})/n).
\end{eqnarray*}
Using first the Cauchy inequality, and then extending the range of
summation over $g$, we derive that
\begin{eqnarray*}
\lefteqn{\left|\sum_{\substack{g=1\\
\gcd(g,n) =1}}^n  \sum_{\substack{1\le h\le n\\ \gcd(h,n)=1}}
\e(a(g^n - g h^{s_n})/n)\right|^2} \\
&& \qquad\le \varphi(n) \sum_{g=1}^n \left| \sum_{\substack{1\le
h\le n\\ \gcd(h,n)=1}} \e(ag h^{s_n}/n)\right|^2 =  \varphi(n)\,n
M_a(n,s_n),
\end{eqnarray*}
where
$$
M_a(n,s) = \#\{(x,y) ~:~  ax^s \equiv ay^s \pmod n,~ x,~y \in
(\Z/n\Z)^*\}.
$$
Now, clearly $M_a(n,s) = \varphi(n) L_a(n,s)$, where
$$
L_a(n,s) = \#\{x ~:~ ax^s \equiv a \pmod n, ~x  \in (\Z/n\Z)^*\}.
$$
Therefore,
\begin{equation}
\label{eq:last2} |\sigma|\le  \sum_{n\in \cN} \sqrt{ n
L_a(n,s_n)}.
\end{equation}

Since $n\in \cN$, there exists a prime $q\in \cI$ such that
  $q\mid d_n$ but $q^2 \nmid n$. Let $\alpha \ge 1$ be the largest
power of $q$ dividing $\lambda(n)$. Then there exists prime
  $p\mid n$ such that $q^{\alpha} \mid p-1$. It is also clear that
$q^\alpha \nmid s_n$. This immediately shows that $\gcd(s_n,
p-1)\mid(p-1)/q$. Since, by the Chinese Remainder Theorem,
$L_a(n,s)$ is a multiplicative function with respect to $n$ (and
since $p>q
>  y$ we also have both $\gcd(n/p, p) =1$ and $\gcd(a, p) =1$), we
derive that
\begin{eqnarray*}
L_a(n,s_n) & = & L_a(n/p,s_n) L_a(p,s_n) \le \varphi(n/p)
L_a(p,s_n)
= \varphi(n/p)  L_1(p,s_n) \\
& = & \varphi(n/p)  \gcd(s_n, p-1) \le \varphi(n/p) (p-1)/q =
\varphi(n)/q \le n/y.
\end{eqnarray*}
Now the relation~\eqref{eq:last2} immediately shows that $\sigma
\ll N^2 y^{-1/2}$, which together with~\eqref{eq:W and sigma}
concludes the proof.
\end{proof}

\section{Open Questions}

Clearly, the range over $a$ in Theorems~\ref{thm:Sum S},
\ref{thm:Sum D} and~\ref{thm: D'} can easily be extended. However,
we do not see how to improve the corresponding bounds, even at the
cost of reducing the range of $a$. Neither can we see any
approaches toward estimating the single exponential sums
\begin{eqnarray*}
T_g(a;N) &=& \sum_{\substack{n=1\\
n~{\mathrm{composite}}}}^N \e(a f_g(n)),\\
{\widetilde T}_g(a;N) &=& \sum_{\substack{n=1\\
n~{\mathrm{composite}}}}^N \e(a {\widetilde f}_g(n)),
\end{eqnarray*}
and we would like to leave these as open problems.

\end{document}